\documentclass[11pt]{amsart}

\usepackage{epsf,epsfig,amsfonts,amsgen,amsmath,amstext,amsbsy,amsopn,amsthm,cases,listings,color}

\usepackage{ebezier,eepic}
\usepackage{mathrsfs}
\usepackage{enumitem}
\usepackage{booktabs}
\usepackage{mathtools}
\usepackage{amssymb}
\usepackage{mathrsfs}

\allowdisplaybreaks[1]

\newtheorem{definition}{Definition} [section]
\newtheorem{theorem}[definition]{Theorem}

\newtheorem{conjecture}[definition]{Conjecture}

\newtheorem{problem}[definition]{Problem}

\newenvironment{pf}{\noindent{\bf Proof}}

\def\qed{\hfill \rule{4pt}{7pt}}

\renewcommand{\eps}{\varepsilon}
\newcommand{\cC}{\mathcal{C}}
\newcommand{\cE}{\mathcal{E}}
\newcommand{\cH}{\mathcal{H}}
\newcommand{\cM}{\mathcal{M}}

\newcommand{\cK}{\mathcal{K}}
\newcommand{\cP}{\mathcal{P}}
\newcommand{\cT}{\mathcal{T}}
\newcommand{\ccW}{\mathscr{W}}

\newcommand{\bP}{\mathbb{P}}
\newcommand{\bE}{\mathbb{E}}

\usepackage[english]{babel}

\usepackage[letterpaper,top=2cm,bottom=2cm,left=3cm,right=3cm,marginparwidth=1.75cm]{geometry}

\usepackage{amsmath}
\usepackage{graphicx}
\usepackage[colorlinks=true, allcolors=blue]{hyperref}

\title{Ramsey theory constructions from hypergraph matchings}

\author{Felix Joos}
\address{Institut f\"ur Informatik, Universit\"at Heidelberg, Germany.
Research supported by the Deutsche Forschungsgemeinschaft (DFG, German Research Foundation) -- 428212407}
\email{joos@informatik.uni-heidelberg.de}

\author{Dhruv Mubayi}
\address{Department of Mathematics, Statistics, and Computer Science, University of Illinois, Chicago, IL, 60607 USA. 
Research partially supported by NSF grants DMS-1763317, 1952767, 2153576 and a Humboldt Research Award}
\email{mubayi@uic.edu}

\begin{document}

\begin{abstract}
 We give asymptotically optimal constructions in generalized Ramsey theory using results about conflict-free hypergraph matchings. 
 For example, we present an edge-coloring of $K_{n,n}$ with $2n/3 + o(n)$ colors such that each $4$-cycle receives at least three colors on its edges. 
 This answers a question of Axenovich, F\"uredi and the second author (On generalized Ramsey theory: the bipartite case,
 J.\ Combin.\ Theory Ser\ B 79 (2000), 66--86). 
 We also exhibit an edge-coloring of $K_n$ with $5n/6+o(n)$ colors that assigns each copy of $K_4$ at least five colors. This gives an alternative very short solution to an old question of Erd\H os and Gy\'arf\'as that was recently answered by Bennett, Cushman, Dudek, and  Pra\l{}at by analyzing a colored modification of the triangle removal process.  
\end{abstract}

\maketitle

\section{Introduction}

Given graphs $G$ and $H$, an $(H,q)$-coloring of~$G$ is an edge-coloring of~$G$ such that every copy of~$H$ in~$G$ receives at least distinct~$q$ colors. 
Let $r(G, H, q)$ be the
minimum number of colors in an $(H,q)$-coloring of $G$.  
 Classical Ramsey numbers
 for multicolorings are  the special case $G=K_n, H=K_p, q=2$.
Initiated by Erd\H os and Shelah \cite{Erd:75, Erd:81}, and subsequently developed by Erd\H os and Gy\'arf\'as~\cite{EG:97} and by Axenovich, F\"uredi and the second author~\cite{AFM:00}, this
generalization of Ramsey numbers has given rise to many interesting problems that have been studied over the years~\cite{A, CFLS:15, Mub:98, Mub:04}. 

One of the main  problems posed in~\cite{AFM:00} was to determine the asymptotics of the smallest open bipartite case when $q>2$, namely, $r(K_{n,n}, C_4, 3)$. In this paper we answer this question by proving 
\begin{equation} \label{eqn:knn}
    r(K_{n,n}, C_4, 3) = 2n/3+o(n).
\end{equation} Since the work of~\cite{AFM:00}, it was evident why proving (\ref{eqn:knn}) would be challenging. Indeed, the lower bound argument shows that the required coloring must arise from an asymptotically optimal  resolvable bipartite partial Steiner quadruple system, where the union of any two color classes has girth at least five (given an appropriate notion of girth). The girth condition between color classes is the most difficult part of this construction.
Using very recent work of  Glock, Kim, K\"uhn, Lichev, and the first author~\cite{GJKKL:22} on conflict-free hypergraph matchings we are able to provide such a construction. 

 Our method also yields 
 \begin{align}\label{eq:56n}
r(K_n, K_4, 5)= 5n/6 + o(n).     
 \end{align}
This gives a new proof of a very recent result of 
Bennett, Cushman, Dudek, and Pra\l{}at~\cite{BCDP:22} which answered a question of Erd\H os and Gy\'arf\'as~\cite{EG:97}. 
In~\cite{BCDP:22},
the authors analyze a modification of the triangle removal process
where each selected triangle receives two colors (one which is used twice and one which is used once).
In addition, certain vertices are not allowed to be incident to a given set of colors and at no point of the process is a 2-colored 4-cycle  allowed to be created.
Not surprisingly, the analysis of this process is technically involved and challenging.

One of the main messages of this paper is that such an intricate analysis is not needed,
as the main result in~\cite{GJKKL:22} implies results of this type even in this fairly involved and particular setting with colors and further constraints. This is achieved by translating 
 the edge-coloring problem to the problem of obtaining a large matching in a suitable auxiliary hypergraph.
Presenting this translation is one of our main contributions.
We are confident that this will be helpful in many other scenarios as well.
To be clear,
it is still probably helpful to propose a random process that creates the desired outcome with high probability, but the analysis may not be needed as the process may be described as a hypergraph matching process and then the result in~\cite{GJKKL:22} can be applied.

\section{Hypergraph matchings with conflicts and other preliminaries}

In this section we state a simplified version of the conflict-free hypergraph matching theorem from~\cite{GJKKL:22},
which is the main tool that we use in our proofs.
Roughly speaking, this result says that vertex-regular uniform hypergraphs $\cH$ with small codegrees admit almost perfect matchings that avoid specified edge conflicts.
Here, conflicts are sets of disjoint edges, modelled as edges of a hypergraph $\cC$ with vertex set $E(\cH)$, which  could, in principle, appear as a subset of the matching. 
The theorem guarantees that all these subsets of edges can be avoided in an almost perfect matching 
provided the conflicts satisfy some natural (and necessary) conditions.

We refer to a $k$-uniform hypergraph as a $k$-graph and identify the edge set with the hypergraph.
For a not necessarily uniform hypergraph $\cC$ and an integer $k$, we use $\cC^{(k)}$ to denote the subgraph of $\cC$  comprising all edges of size exactly $k$.
For a vertex $e$ of $\cC$, we write $\cC_e$ for the set of all $C\setminus\{e\}$ with $C\in \cC$ and $e\in C$.
For a hypergraph $\cH$ and a vertex $v$ of $\cH$, 
we use $d_\cH(v)$ to denote the \emph{degree} of $v$; 
that is, the number of edges containing $v$.
We use $\Delta(\cH)=\max_{u\in V(\cH)}d_\cH(u)$ to denote the \emph{maximum degree} of $\cH$ and,
similarly, we define the \emph{minimum degree} $\delta(\cH)$ of $\cH$.
For $j\geq 2$, we denote by~$\Delta_j(\cH)$ the maximum number of edges that contain a particular set of $j$ vertices.
We also define $[k]=\{1,\ldots,k\}$.
We use standard asymptotic notation $o(\cdot),O(\cdot), \Theta(\cdot)$ as $n \to \infty$, which we always treat as (a set of) non-negative functions. 

For a hypergraph $\cH$, we refer to $\cC$ as a \emph{conflict system} for $\cH$ if $\cC$ is hypergraph with vertex set $E(\cH)$.
We say that a set $E$ of edges of $\cH$ is \emph{$\cC$-free} if no $C\in \cC$ is a subset of $E$.
For integers~$d\geq 1, \ell\geq 3$ and~$\eps\in(0,1)$, 
we say that~$\cC$ is \emph{$(d,\ell,\eps)$-bounded} if the following holds.
\begin{enumerate}[label=\textup{(C\arabic*)}]
		\item\label{item: conflict size} $3\leq|C|\leq \ell$ for all~$C\in \cC$;
		\item\label{item: conflict degree} $\Delta(\cC^{(j)})\leq \ell d^{j-1}$ for all $3\leq j \leq \ell$;
		\item\label{item: conflict codegrees} $\Delta_{j'}(\cC^{(j)})\leq d^{j-j'-\eps}$ for all~$3\leq j \leq \ell$ and~$2\leq j' \leq j-1$.
\end{enumerate}
This means $\cC$ is $(d,\ell,\eps$)-bounded if the conflicts have size at most $\ell$
and that for all uniform subgraphs of~$\cC$ the degrees and codegrees are not too large.

Suppose now that $\cC$ is $(d,\ell,\eps$)-bounded. 
The following concept of weight functions can be used to ensure that the provided almost perfect matching admits quasirandom properties.
A \emph{test function} for~$\cH$ is a function~$w\colon \binom{\cH}{j}\rightarrow[0,\ell]$ where~$j\in  \mathbb{N}$ such that~$w(E)=0$ whenever~$E\in\binom{\cH}{j}$ is not a matching.
We refer to~$j$ as the uniformity of~$w$ and we say that~$w$ is $j$-uniform.
In general, for a function~$w\colon A\rightarrow\mathbb{R}$ and a finite set~$X\subset A$, we define~$w(X):=\sum_{x\in X}w(x)$.
If~$w$ is a~$j$-uniform test function, we also use~$w$ to denote the extension of~$w$ to arbitrary subsets of~$\cH$ such that for all~$E\subset\cH$, we have~{$w(E)=w(\binom{E}{j})$}.
For~$j,d\in \mathbb{N}$, $\eps>0$ and a conflict system~$\cC$ for~$\cH$, we say that a~$j$-uniform test function~$w$ for~$\cH$ is~\emph{$(d,\eps,\cC)$-trackable} if the following holds. 
\begin{enumerate}[label=\textup{(W\arabic*)}]
		\item\label{item: trackable size} $w(\cH)\geq d^{j+\eps}$;
		\item\label{item: trackable degrees} $w(\{ E\in\binom{\cH}{j} \colon E\supseteq E' \} )\leq w(\cH)/d^{j'+\eps}$ for all~$j'\in[j-1]$ and~$E'\in\binom{\cH}{j'}$;
		\item\label{item: trackable neighborhood} $|(\cC_e)^{(j')}\cap (\cC_f)^{(j')}|\leq d^{j'-\eps}$ for all~$e,f\in \cH$ with~$w(\{E\in\binom{\cH}{j}\colon e,f\in E\})>0$ 
		and all~$j'\in[\ell-1]$;
		\item\label{item: trackable no conflicts} $w(E)=0$ for all~$E\in\binom{\cH}{j}$ that are not~$\cC$-free.
\end{enumerate}

\begin{theorem}[{\cite[Theorem~3.3]{GJKKL:22}}]\label{hypmat}
	For all~$k,\ell\geq 2$, there exists~$\eps_0>0$ such that for all~$\eps\in(0,\eps_0)$, there exists~$d_0$ such that the following holds for all~$d\geq d_0$.
	Suppose~$\cH$ is a $k$-graph on~$n\leq \exp(d^{\eps^3})$ vertices with~$(1-d^{-\eps})d\leq \delta(\cH)\leq\Delta(\cH)\leq d$ and~$\Delta_2(\cH)\leq d^{1-\eps}$ and suppose~$\cC$ is a~$(d,\ell,\eps)$-bounded conflict system for~$\cH$.
	Suppose~$\ccW$ is a set of~$(d,\eps,\cC)$-trackable test functions for~$\cH$ of uniformity at most~$\ell$ with~$|\ccW|\leq \exp(d^{\eps^3})$.
	Then, there exists a $\cC$-free matching~$\cM\subset \cH$ of size at least~$(1-d^{-\eps^3})n/k$ with~$w(\cM)=(1\pm d^{-\eps^3})d^{-j}w(\cH)$
	for all~$j$-uniform~$w\in\ccW$.
\end{theorem}

\section{$r(K_{n,n}, C_4, 3)$}

It is shown in~\cite{AFM:00} that $r(K_{n,n}, C_4, 3)\geq  2n/3 $.
Hence it order to prove~\eqref{eqn:knn}, it suffices to show that there exists $\delta>0$ for which $r(K_{n,n}, C_4, 3)\leq  2n/3+ n^{1-\delta}$ by providing a coloring of~$K_{n,n}$ with at most $2n/3+ n^{1-\delta}$ colors such that any $4$-cycle in $K_{n,n}$ receives at least three colors.
Our construction is obtained in two stages. 
In the first stage, we translate the problem to a question about hypergraph matchings and then apply Theorem~\ref{hypmat}. 
This produces an edge-coloring of all but $O(n^{2-\delta})$ edges of $K_{n,n}$. 
Moreover, the distribution of the uncolored edges admits certain quasirandomness properties.
In the second stage, we color the uncolored edges randomly with a new set of colors and ensure, by the local lemma, that all $4$-cycles receive at least three colors.

We proceed with the first part.
To this end, we need the following notation.
For $k \ge 2$, a cycle of length~$k$ in a hypergraph is a collection of $k$ distinct edges $\{e_1, \ldots, e_k\}$ such that there exist~$k$ distinct vertices $v_1, \ldots, v_k$ with $e_i \cap e_{i+1} = \{v_i\}$ for all $i\in [k]$, where indices are taken modulo~$k$. 
The girth of a hypergraph $\cH$ is the minimum length of a cycle in~$\cH$. 
For disjoint sets $X,Y$, we write $B(X, Y)$ for the complete bipartite graph with parts $X, Y$.

\begin{theorem} \label{stage1}
There exists $\delta>0$ such that for all sufficiently large $n$ in terms of $\delta$, there exists an edge-coloring of a subgraph of $G \subset B(X, Y) \cong K_{n,n}$ with at most $2n/3$ colors and the following properties:
\begin{enumerate}[label=\rm (\Roman*)]
    \item Every color class consists of vertex-disjoint 3-edge stars. 
    \item Given any two colors $i,j$, 
    the 4-graph with vertex set $X\cup Y$ where each edge is formed by the vertex set of a star in color $i$ or $j$ has girth at least 5.
    \item The graph $L = B(X, Y)-E(G)$ has maximum degree at most $n^{1-\delta}$. 
    \item For each $(x,y) \in X \times Y$, the number of  $x'y'\in E(L)$ with $x'\in X\setminus \{x\},y'\in Y\setminus \{y\}$ 
    such that $xy',yx'$ receive the same color is at most $n^{1-\delta}$.
\end{enumerate}
\end{theorem}
\proof  
Our aim is to apply Theorem~\ref{hypmat}.
To this end, we define the following 10-uniform hypergraph~$\cH$ as follows.
Let $X,Y$ be two disjoint vertex sets of size $n$.
The vertex set of~$\cH$ is $U \cup V$, where 
$U = {X \cup Y \choose 2}$ and
$V = \bigcup _{i\in [2n/3]} V_i$ where each $V_i$ is a copy of $V(B(X, Y))=X \cup Y$. 
Let $\cT$ be the collection of 4-sets in $X\cup Y$ which have exactly one or three vertices in $X$. 
Given $e \in \cT$ and $i \in [2n/3]$, 
let $e_i= {e \choose 2} \cup e'_i$ where $e'_i$ is the copy of $e$ in $V_i$. 
Thus $e_i$ has six vertices in $U$ and four vertices in $V_i$.
Finally, let 
$$E(\cH) = \{e_i: e \in \cT, i \in [2n/3]\}.$$
We write $\cK$ for the set of all copies of $K_4$ in the complete graph with vertex set $X \cup Y$ which have exactly one or three vertices in $X$.
Observe that there is a natural bijection between the edge set of $\cH$ and the set of tuples $(K,i)$ where $K\in \cK$ and $i \in [2n/3]$,
namely, by mapping $e_i\in \cH$ with $e\in \cT,i \in [2n/3]$ to $(K,i)$ where $K$ is the complete graph on $e$ of size~$4$.
Moreover,  there is also a bijection between $\cK$ and all $3$-edge stars in $B(X,Y)$ by simply deleting three edges from $K\in \cK$ or adding the three missing edges to a star.
Hence we refer to an edge $e$ in $\cH$ also as a $K_4$ or star with color $i$ and simply write $e=(K,i)$ for $K\in \cK$ and $i\in[2n/3]$.
Using this bijection
each matching in $\cH$ corresponds to an edge-coloring of a subgraph of~$B(X,Y)$ with at most $2n/3$ colors satisfying (I).
In the following, we aim to find a particular matching which also satisfies (II)--(IV).

 We claim that $\cH$ is essentially $d$-regular with $d=2n^3/3$.
 Indeed, fix $u \in U$. 
 If $u\in {X \choose 2}$, then
 in order to pick an edge containing~$u$, 
 we pick a vertex $x \in X \setminus u$, 
 a vertex $y \in Y$ and $i \in [2n/3]$. 
 This gives $(|X|-2)|Y|(2n/3) = d- O(n^2)$ edges. 
 The same argument works for  $u \in {Y \choose 2}$. 
 Now suppose that  $u =\{x,y\}$ with $x\in X,y\in Y$.
 In this case we pick $Z\in \{X,Y\}$ and then distinct $z, z' \in Z\setminus u$ and $i \in [2n/3]$ to obtain $2{n-1 \choose 2}(2n/3) = d- O(n^2)$ edges.  
 Thus in any case $$d - O(n^2) \leq d_\cH(u) \leq  d. $$
 Now fix $v \in V_i$. 
 Edges containing $v$ in $\cH$ arise only from copies of $K_4$ containing the copy of~$v$ in $B(X,Y)$ and the number of such copies is 
 ${n-1 \choose 2}n +{n \choose 3} = d - O(n^2).$
 Therefore, 
 $$d\left(1-\frac{1}{n^{1/2}}\right) \le \delta(\cH) \le \Delta(\cH) \le d.$$ 
 It is easy to prove that $\Delta_2(\cH) \le n^2$.
 
 We next define a 4-graph $\cC$ which is a conflict system for $\cH$. 
 As required, we set $V(\cC)=E(\cH)$. 
Edges of $\cC$ arise from $4$-cycles in $B(X,Y)$ comprising two monochromatic matchings of size~2 (hence resulting in a $4$-cycle with exactly two colors). 
 Such 2-colored 4-cycles arise from four stars in $B(X,Y)$, two in one color and two in another color, that form a (linear)
 4-uniform 4-cycle. 
 More precisely, given four distinct vertices $x,x' \in X$ and $y,y' \in Y$, two distinct colors $i,j\in [2n/3]$, and $e=(K_{xy},i),e'=(K_{x'y'},i),f=(K_{xy'},j),f'=(K_{x'y},j)$
 with $a,b \in V(K_{ab})$ for $a\in \{x,x'\}, b \in \{y,y'\}$, we have $\{e, e', f, f'\} \in E(\cC)$ whenever $\{e, e', f, f'\}$ is a matching in~$\cH$.
 
 We next claim that $\Delta(\cC) = O(d^3)$. 
 Indeed, given a vertex $e_i=(K,i) \in V(\cC)$,
 all edges containing~$e_i$ are of the form $\{(K,i),(K',i),(K_1,j),(K_2,j)\}$ for some $K',K_1,K_2\in \cK$ and $j\in [2n/3]\setminus \{i\}$.
 There are $O(n)$ choices for $j$, $O(n^4)$ choices for $K'$, and $O(n^4)$ choices for $V(K_1)\cup V(K_2)$. 
 Hence $\Delta(\cC)=O(n^9)=O(d^3)$.
 
 A similar calculation also shows that there are $O(n^5)$ edges in $\cC$ containing any fixed $(K,i),(K',i)\in V(\cC)$.
 Similarly, if we fix $(K,i), (K_1,j)$, there are $O(n^3)$ choices for $(K',i)$ and then  $O(n^2)$ choices for $(K_2,j)$. 
 Consequently, $\Delta_2(\cC) = O(n^5) < d^{2-1/4}$ with room to spare. 
 Given $(K,i),(K',i),(K_1,j)$ as above leaves $O(n^2)$ choices for $(K_2,j)$.
 Thus $\Delta_3(\cC)<d^{1-1/4}$  and~$\cC$ is a $(d, O(1), \varepsilon)$-bounded conflict system for~$\cH$ for all $\varepsilon \in (0, 1/4)$.
 
 We remark that one can now apply Theorem~\ref{hypmat} to $\cH$ to obtain an almost perfect conflict-free matching $M$, and this translates to an edge-coloring of a subgraph $G$ of most edges of $B(X,Y)$  with every $4$-cycle receiving at least three colors (each color class is a star forest with stars of three edges);
 in particular, Theorem~\ref{stage1} (I) and (II) hold. 
 However, to complete this to a coloring of all of $B(X,Y)$, we need to ensure that the coloring of $G$ has further properties (properties (III) and (IV)). 
 We achieve this by considering carefully chosen test functions.
 
To this end, for each  $v \in X \cup Y$, let $S_v\subset U$ be the set of edges in $B(X,Y)$ incident to~$v$.
 Let $w_v:E(\cH)\to [0,3]$ be the weight function that assigns every edge of $\cH$ the size of its intersection with $S_v$.
 Note that $w_v(M)$ counts the number of edges in $S_v$ that belong to a star in any collection of stars that stems from some matching~$M$ in~$\cH$.
 Moreover,  $w_v(\cH)=\sum_{e \in S_v}d_\cH(e)=nd-O(n^3)$.
 Observe that~\ref{item: trackable degrees}--\ref{item: trackable no conflicts} are trivially satisfied, because~$w_v$ is $1$-uniform.
 Hence $w_v$ is a $(d, \varepsilon, \cC)$-trackable 1-uniform test function for all $\varepsilon\in (0,1/4)$ and $v \in V(B(X,Y))$.

We could now apply Theorem~\ref{hypmat} to also obtain (III). 
Indeed, for suitable $\varepsilon$, and sufficiently  large~$n$,
 Theorem~\ref{hypmat} yields a $\cC$-free matching~$M\subset \cH$ such that for each $v \in X \cup Y$, 
 \begin{align*}
 w_v(M)>(1 - d^{-\eps^3})d^{-1}w_v(\cH)
   >(1 - n^{-\delta})n
 \end{align*}
 for some small enough $\delta>0$.
 Therefore, for every vertex $v\in X\cup Y$, there are at most $n^{1-\delta}$ edges in $B(X,Y)$ incident to $v$ that do not belong to a star selected by $M$.

Next we extend our selection of test functions such that Theorem~\ref{hypmat} also yields (IV).
For each $(x,y) \in X \times Y$, we proceed as follows.
For each $j_x,j_y\in \{1,3\}$, let
\begin{align*}
\cP_{j_x,j_y}
&=\{\{(K_x,i),(K_y,i)\}\colon i\in [2n/3], V(K_x)\cap V(K_y)=\emptyset,\\ 
&x\in V(K_x), |V(K_x)\cap X|=j_x,y\in V(K_y), |V(K_y)\cap Y|=j_y\}.
\end{align*}

Clearly, 
$$|\cP_{j_x,j_y}|=\frac{j_xn^3}{6} \cdot \frac{j_yn^3}{6} \cdot \frac{2n}{3}  \pm O(n^6)= \frac{j_xj_yn^7}{54} \pm O(n^6) > d^{2+1/4}.$$
We define an indicator weight function $w_{x,y,j_x,j_y}$ for the pairs in $\cP_{j_x,j_y}$.
Assume for now  that these test functions are $(d, \varepsilon, \cC)$-trackable for all $\varepsilon \in (0, 1/4)$ (note that each pair in~$\cP_{j_x,j_y}$ corresponds to a matching of size $2$ in $\cH$).
Adding all these weight functions (for all $(x,y)\in X\times Y$ and $j_x,j_y\in \{1,3\}$) to the set of weight functions which we give to Theorem~\ref{hypmat}, we obtain a matching $M$ such that 
\begin{equation} \label{eqP}\left|{M \choose 2} \cap \cP_{j_x,j_y}\right| = w_{x,y,j_x,j_y}(M) \leq (1+d^{-\varepsilon^3})d^{-2}|\cP_{j_x,j_y}| \le (1+n^{-2\delta})\frac{j_xj_yn}{24}\end{equation} for each $x,y,j_x,j_y$.

Now we verify that $w_{x,y,j_x,j_y}$ is indeed $(d, \varepsilon, \cC)$-trackable. 
We already checked \ref{item: trackable size}.
To see~\ref{item: trackable degrees}, we fix one $(K_x,i)$
and observe that there are $O(n^3) <n^7/d^{1+1/4}$ choices for $(K_y,i)$ (and similarly when $x,y$ are swapped), which shows~\ref{item: trackable degrees}.
To see~\ref{item: trackable neighborhood}, fix $e=(K_1,i),f=(K_2,j)$
with $w_{x,y,j_x,j_y}(\{e,f\})>0$.
Hence $i=j$ and, say, $e=(K_x,i),f=(K_y,i)$ as in the definition of~$\cP_{j_x,j_y}$.
To choose the three edges of $\cH$ that form a conflict with both $e$ and $f$,
we can choose at most six further vertices outside $V(K_x)\cup V(K_y)$ and one color in $[2n/3]\setminus \{i\}$
resulting in $O(n^7)<d^{3-1/4}$ choices, which shows~\ref{item: trackable neighborhood}.
Property~\ref{item: trackable no conflicts} is vacuously true.

Now we define a set of triples $\cT_{j_x,j_y}$ that extend the pairs in $\cP_{j_x,j_y}$.
Given $\{(K_x,i),(K_y,i)\}\in \cP_{j_x,j_y}$, we add the triple $\{(K_x,i),(K_y,i),(K,j)\}$ to $\cT_{j_x,j_y}$ whenever $j\in [2n/3]\setminus \{i\}$ and $K$ contains exactly one vertex in $V(K_{x})\setminus X$ and exactly one vertex in $V(K_{y})\setminus Y$.
Consequently, we have $$|\cT_{j_x,j_y}|=(4-j_x)(4-j_y)(d\pm O(n^2))|\cP_{j_x,j_y}|$$
and all triples in~$\cT_{j_x,j_y}$ are matchings of size $3$ in $\cH$.
In the same vein as above, we define an indicator weight function $w_{x,y,j_x,j_y}'$ on the triples of~$E(\cH)$.
Assume for now that these weight functions are $(d, \varepsilon, \cC)$-trackable for all $\varepsilon \in (0, 1/4)$.
 Theorem~\ref{hypmat} gives rise to a matching $M$ such that 
\begin{equation} \label{eqT}\left|{M \choose 3} \cap \cT_{j_x,j_y}\right|\geq (1-n^{-2\delta})(4-j_x)(4-j_y) \cdot \frac{j_xj_yn}{24}\end{equation}
for each $x,y,j_x,j_y$.

The crucial observation is that each $x'y'$ from (IV) arises from an edge of $\cH$ in some triple of~$\cT_{j_x,j_y}$ that contains some pair from 
$\cP_{j_x,j_y}$. Therefore the number of $x'y'$ as in (IV) is at most 
$$    \sum_{j_x,j_y\in \{1,3\}} (4-j_x)(4-j_y)\left|{M\choose 2} \cap \cP_{j_x,j_y}\right|- \left|{M\choose 3} \cap \cT_{j_x,j_y}\right|.$$
By (\ref{eqP}) and (\ref{eqT}), we see that Theorem~\ref{hypmat} yields a matching M such that the above quantity at most $n^{1-\delta}$, which proves~(IV).

It remains to check that $w_{x,y,j_x,j_y}'$ is $(d, \varepsilon, \cC)$-trackable for all $\varepsilon \in (0, 1/4)$.
Observe that $w_{x,y,j_x,j_y}'(\cH)=\Theta(n^{10})> d^{3+1/4}$ and hence \ref{item: trackable size} holds.
To see~\ref{item: trackable degrees}, we observe that $w_{x,y,j_x,j_y}'(\{ E\in\binom{\cH}{3} \colon E\supseteq e \} )= O(n^6) <w_{x,y,j_x,j_y}'(\cH)/d^{1 +1/4}$
for all $e\in \cH$
and also $w_{x,y,j_x,j_y}'(\{ E\in\binom{\cH}{j} \colon E\supseteq E' \} )= O(n^3)<w_{x,y,j_x,j_y}'(\cH)/d^{2 +1/4}$
for all $E'\in \binom{\cH}{2}$.
To see~\ref{item: trackable neighborhood}, 
fix $e=(K_1,i),f=(K_2,j)$
with $w_{x,y,j_x,j_y}(E)>0$ for some $\{e,f\}\subset E \in \binom{\cH}{3}$.
If $i=j$ as above for $\cP_{j_x,j_y}$, then we again conclude as above the desired bound.
Hence assume that $i\neq j$.
There at most $O(n^8)<d^{3-1/4}$ conflicts containing $e$ which involve only edges of~$\cH$ with colours $i,j$.
This is clearly an upper bound for $|(\cC_e)^{(3)}\cap (\cC_f)^{(3)}|$ and thus implies~\ref{item: trackable neighborhood}.
Property~\ref{item: trackable no conflicts} is vacuously true, which completes the proof.
\qed

Next we finish the proof of \eqref{eqn:knn}.
First we apply Theorem~\ref{stage1} (with $2\delta$ playing the role of $\delta$) and obtain a coloring as stated in Theorem~\ref{stage1}.
In particular, $\Delta(L)\leq n^{1-2\delta}$ by (III).
We now show how to color the edges of $L$ with a set $P$ of $k=n^{1-\delta}$ new colors so that no $4$-cycle in the final coloring has fewer than three colors on its edges. 
We color each edge of $L$ with a color from~$P$ with equal probability $1/k$, 
independently of all other edges. 
We apply the symmetric form of the local lemma to show that there is no 2-colored $4$-cycle and to this end we now define three types of bad events. 

For any pair $e,f$ of adjacent edges in $L$ and $i\in P$,
we let $A_{e,f,i}$ be the event that both~$e$ and~$f$ receive color~$i$.
Clearly, $\bP[A_{e,f,i}]=k^{-2}$.
For $4$-cycle $D$ in~$L$,
we define $B_D$ to be the event that~$D$ receives exactly two colors (and this coloring is proper).
Hence $\bP[B_{D}]=(k(k-1))^{-1}\leq 2k^{-2}$.
For a $4$-cycle $D=xyx'y'$ in $K_{n,n}$ and $i\in P$,
where $xy,x'y'\in E(L)$ and $xy',yx'$ belong to two stars colored alike by Theorem~\ref{stage1},
we let $C_{D,i}$ be the event that both $xy$, $x'y'$ receive color $i$.
Clearly, $\bP[C_{D,i}]=k^{-2}$.
Let $\cE$ be the collection of all these defined events.

We say two events as above are edge-disjoint if the edges in $L$ that can be associated with them are distinct.
Fix an event $E\in \cE$ as above.
There are at most $8\Delta(L)\cdot  k \leq 8k^2 n^{-\delta}$ events~$A_{e,f,i}$ that are not edge-disjoint from $E$;
there are at most $4(\Delta(L))^2\leq k^2 n^{-\delta}$ events~$B_{D}$ that are not edge-disjoint from $E$; and
there are at most $4n^{1-2\delta}k\leq 4k^2 n^{-\delta}$ events $C_{D,i}$ that are not edge-disjoint from $E$ by (IV).
Consequently,
for any event $E\in \cE$, there is a set~$\cE_E$ of $13k^2 n^{-\delta}$ events from $\cE$ such that
$E$ is independent of any collection of events in $\cE\setminus \cE_E$.
Since $13k^2 n^{-\delta} \cdot (2k^{-2})\leq 1/4$,
by the local lemma, there is a coloring of the edges of $L$ such that none of the events in $\cE$ hold.
Then given such a coloring, it is easy to see that every $4$-cycle in $K_{n,n}$ receives at least three distinct colors.
\qed

\section{$r(K_n, K_4, 5)$}

In this section we prove~\eqref{eq:56n}.
Recently,
\eqref{eq:56n} was proven in~\cite{BCDP:22}
by analyzing an elaborated extension of the triangle removal process.
Here we encode (the result of) this random process as a hypergraph matching problem (with conflicts)
in a similar fashion as in the previous section.
We expect that this methodology may be fruitful for many other applications as well.

It will be convenient to use the following concentration inequality due to McDiarmid.

\begin{theorem}[McDiarmid's inequality, see~\cite{mcdiarmid:89}] \label{thm:McDiarmid}
Suppose $X_1,\dots,X_m$ are independent random variables.
Suppose $X$ is a real-valued random variable determined by $X_1,\dots,X_m$ such that changing the outcome of $X_i$ changes $X$ by at most $b_i$ for all $i\in [m]$.
Then, for all $t>0$, we have 
$$\bP[|X-\bE[X]|\ge t]\le 2 \exp\left({-\frac{2t^2}{\sum_{i\in [m]} b_i^2}}\right).$$
\end{theorem}

We now state and prove the main technical statement which gives a partial coloring of $K_n$.

\begin{theorem}\label{thm:56n}
    There exists $\delta>0$ such that for all sufficiently large $n$ in terms of $\delta$,
    there exists an edge-coloring of a subgraph $F\subset K_n$ with at most $5n/6 + n^{1-\delta}$ colors and the following properties:
    \begin{enumerate}[label=\rm (\Roman*)]
    \item Every color class consists of vertex-disjoint edges and 2-edge paths. 
     \item For all triangles $xyz$ in $F$ where $xy,yz$ receive the same color and $xz$ is colored $i$,
    the vertex $y$ is an isolated vertex and $xz$ forms a component in color class $i$.
    \item Every $4$-cycle in $F$  receives at least three distinct colors.
    \item The graph $L=K_n-E(F)$ has maximum degree at most $n^{1-\delta}$. 
    \item For each $xy \in E(K_n)$, the number of  $x'y' \in E(L)$ with $\{x,y\} \cap \{x', y'\}=\emptyset$ for which $xx'$ and $yy'$ receive the same color in $F$ is at most $n^{1-\delta}$.
\end{enumerate}
\end{theorem}

\proof
Let $\delta>0$ be sufficiently small and let $n$ be sufficiently large in terms of $\delta$.
Let $\rho=n^{-\delta}$ and $k=(1+\rho)5n/6$.
Let $G=K_n$. We first construct a random (vertex) set $V$ as follows.
Let $V'=\bigcup_{i\in [k]}V_i'$ where $V_1',\ldots,V_k'$ are copies of $V(G)$.
Delete each vertex in $\bigcup_{i\in [k]}V_i'$ independently with probability $p=\rho/(1+\rho)$.
Denote by $V_i$ the vertices in $V_i'$ which remain.
Let $V=\bigcup_{i\in[k]}V_i$.

Next, we construct an $8$-graph $\cH$ with vertex set $E(G)\cup V$ as follows.
For $i\in [k]$ and $v\in V(G)$, we denote by $v_i$ the copy of $v$ in $V_i$ (if it exists).
For each triangle $uvw$ in $G$ and distinct $i,j\in [k]$,
we add the edge
\begin{align*}
    \{uv,uw,vw, u_i,v_i,w_i,v_j,w_j\}
\end{align*}
to $\cH$ if $u_i,v_i,w_i,v_j,w_j\in V$ and $u_j\notin V$. 

Next we consider triangles where one vertex has a label.
Clearly, each unlabelled triangle gives rise to three different such labelled triangles.
Let $\cK$ be the set of all such labelled triangles that arise from triangles in $G$.
Observe that there is an injection from the edge set of $\cH$ to the set of tuples $(K,i,j)$ with $K\in \cK$ and distinct $i,j \in [k]$, 
where we think of color $i$ appearing on two edges in~$K$ and color $j$ appearing on the edge not incident to the labelled vertex of~$K$.
Hence we refer to an edge~$e$ in~$\cH$ also as a triangle (if the colors and the label do not matter at that point) and simply write $e=(K,i,j)$ for $K\in \cK$ and distinct $i,j \in[k]$.

Note that any matching $M \subset \cH$ corresponds to a collection of edge-disjoint triangles in~$G$ where each triangle $uvw$ is of the form where $vw$ is assigned color $j$ and $uv$ and $uw$ are both assigned color $i$. 
Moreover, vertex $u$ is not incident to any edge in color $j$ in any other triangle (due to the condition $u_j \not\in V$ above). 
This yields a  $k$-coloring of the edges of $G$ that lie within these triangles that satisfies properties (I) and (II). Our goal is to find a particular $M$ such that the corresponding edge-coloring of $G$ also satisfies (III)--(V). 
Property (III) will follow by defining a particular conflict system $\cC$ for $\cH$ and finding a $\cC$-free matching $M$, and (IV) and (V) will follow by considering certain trackable test functions.

Note that $\Delta_2(\cH) = O(n^2)$ (independent of our random choices).
Next, we compute the expected degrees of the vertices in $\cH$.
Note that $1-p=(1+\rho)^{-1}$.
We obtain for $uv\in E(G)$ 
\begin{align*}
    \bE[d_{\cH}(uv)]
    &=(n-2) \cdot k(k-1) \cdot 3 \cdot (1-p)^5p \\
    &= \frac{5}{2}n^2k (1-p)^4p \pm O(n^2).
\end{align*}
and for $u_i\in V_i'$
\begin{align*}
    \bE[d_{\cH}(u_i)\mid u_i\in V_i]
    &=\left(\binom{n-1}{2} + (n-1)(n-2) + (n-1)(n-2)\right) (k-1)(1-p)^4p\\
    &= \frac{5}{2}n^2k (1-p)^4p.
\end{align*}
We claim that by Theorem~\ref{thm:McDiarmid}
\begin{align*}
    \frac{5}{2}n^2k (1-p)^4p -O(n^{8/3}) \leq \delta(\cH)\leq \Delta (\cH) \leq \frac{5}{2}n^2k (1-p)^4p + O(n^{8/3}),
\end{align*}
holds with high probability.
Indeed, fix some $uv\in E(G)$ and for $w\in V'$, define $b_w=n^2$ if~$w$ is a copy of $u$ or $v$ and $b_w=n$ otherwise.
Then $\sum_{w\in V'}b_w^2=O(n^5)$, which shows concentration in an interval of length $O(n^{5/2})$ using Theorem~\ref{thm:McDiarmid}.
Similarly, fix $u_i\in V_i'$ and for $w\in V'\setminus\{u_i\}$, define $b_w=n^2$ if $w$ is a copy of $u$ or $w\in V_i'$ and $b_w=n$ otherwise;
again the same conclusion holds.

Hence $\cH$ is essentially $d$-regular for some $d=\Theta(n^{3-\delta})$. 
In addition, it is easy to check using McDiarmid's inequality that given $u_i,v_i\in V$,
there are $O(pn^2) =O(n^{2-\delta})=O(d/n)$ edges containing both $u_i,v_i$ with high probability 
(indeed, for $w\in V'\setminus \{u_i,v_i\}$, define $b_w=n$ if $w$ is a copy of $u$ or $v$ or $w\in V_i'$ and $b_w=1$ otherwise).
We need that one further quantity is close to its expected value, namely the size of $\cP_{j_x,j_y}$ which will be defined  later in the proof. For now we assume that  that it is close to its expected value with high probability.
From now on, we fix one choice for $\cH$ that satisfies these properties and refer to this deterministic $8$-graph again by $\cH$.
Moreover, we set $d=\Delta(\cH)$ (and hence $\cH$ is essentially $d$-regular and $d=\Theta(n^{3-\delta})$).

We next define a 4-graph $\cC$ with vertex set $E(\cH)$ which is a conflict system for $\cH$. 
The edges of $\cC$ arise from $4$-cycles in $G$ comprising two monochromatic matchings of size~2 (hence resulting in a $4$-cycle in $G$ with exactly two colors). 
 We only consider such 2-colored $4$-cycles that arise from four triangles in $G$, two of which contain one color and the other two contain another color. 
 More precisely, given four distinct vertices $w,x,y,z \in V(G)$ and two distinct colors $i,j\in [k]$ with $e_{uv}= (K_{uv}, \alpha_{uv}, \beta_{uv})$   for all $uv \in \{wx,xy,yz,wz\}$, we have $\{e_{wx}, e_{xy}, e_{yz},  e_{wz}\} \in \cC$ whenever  
 
 \begin{itemize}
    \item color $i$ appears on the edges $wx$ in $K_{wx}$  and $yz$ in $K_{yz}$,
    \item color $j$ appears on the edges $xy$ in $K_{xy}$  and $wz$ in $K_{wz}$,
    \item $e_{wx}, e_{xy}, e_{yz}, e_{wz}$ form a matching in $\cH$.
 \end{itemize}
 We emphasize that we do not stipulate whether color $i$ appears as a single edge or a 2-edge path within the triangle $K_{wx}$ and similarly for the other colors and triangles. 
 We also do not require~$K_{wx}$ and~$K_{yz}$  (and~$K_{xy}$ and~$K_{wz}$) to be vertex-disjoint.  
 We note that a $\cC$-free matching $M$ in~$\cH$ gives rise to an edge-coloring in which each $4$-cycle receives at least three colors. 
 Indeed, by the definition of $\cH$, no color class has a 3-edge path and no two monochromatic 2-edge paths share two vertices.   The only other possibility for a 2-colored $4$-cycle is if it comprises two monochromatic matchings and this is precluded by the definition of~$\cC$. 

It is easy to see that $\Delta(\cC)= O(d^3)$. 
Indeed, fix some edge in $\cH$ that assigns color $i$ to  $wx \in E(G)$.
We have at most $n^2$ choices to fix two further vertices $y,z$ in $G$.
There are~$O(d/n)$ choices for the edges that contain $y_i,z_i$,
at most~$d$ choices for the edges that contain $xy$ (say, $xy$ receives color $j$) 
and~$O(d/n)$ choices for the edges that contain $z_j,w_j$.
With almost the same arguments we also obtain that $\Delta_2(\cC) = O(d^2/n)<d^{2-1/4}$ and $\Delta_3(\cC)= O(d/n)<d^{1-1/4}$.
Consequently, $\cC$ is a $(d, O(1), \varepsilon)$-bounded conflict system for $\cH$ for all $\varepsilon \in (0, 1/4)$.
 
Next we define a set of $(d, \varepsilon, \cC)$-trackable test functions that enables us to obtain properties~(IV) and~(V).
For each  $v \in V(G)$, let $S_v\subset E(G)$ be the set of $n-1$ edges in $G$ incident to~$v$.
Let $w_v:E(\cH)\to [0,2]$ 
be the weight function that assigns every edge of $\cH$ the size of its intersection with $S_v$.
Note that $w_v(M)$ counts that number of edges in $S_v$ that belong to a triangle in any collection of triangles that stems from some matching $M$ in $\cH$.
Moreover,  $w_v(\cH)=\sum_{e \in S_v}d_\cH(e)=nd-O(n^3)$.
Observe that~\ref{item: trackable degrees}--\ref{item: trackable no conflicts} are trivially satisfied, because~$w_v$ is $1$-uniform.
Hence $w_v$ is a $(d, \varepsilon, \cC)$-trackable 1-uniform test function for all $\varepsilon \in (0, 1/4)$ and $v \in V(G)$.
 
We could now apply Theorem~\ref{hypmat} to also obtain (IV). 
Indeed, for suitable $\eps$, and sufficiently large $n$, Theorem~\ref{hypmat} yields a $\cC$-free matching~$M\subset \cH$ such that for each $v \in V(G)$, 
\begin{align*}
 w_v(M)>(1 - d^{-\eps^3})d^{-1}w_v(\cH)
   >(1 - n^{-\delta})n
\end{align*}
holds where we choose $\delta>0$ small enough in terms of $\eps$.
Therefore, for every vertex $v\in V(G)$, there are at most $n^{1-\delta}$ edges in $G$ incident to $v$ that do not belong to a triangle selected by~$M$. This proves (III) and (IV).

We turn to (V). 
In what follows, if $u \in V(G)$, then the notation $K_u$ for a triangle $K_u$ implies that $K_u$ contains  $u$.
For all distinct $x,y \in V(G)$ and $j_x,j_y\in \{1,2\}$, we define 
\begin{align*}
\cP_{j_x,j_y}
&=\{\{(K_x,\alpha_x, \beta_x),(K_y,\alpha_y, \beta_y)\}\colon 
(K_x,\alpha_x, \beta_x),(K_y,\alpha_y, \beta_y)\in \cH \text{ are disjoint},\\ 
&\text{some color $i$ is incident to $z$ exactly $j_z$ times in $K_z$ for each } z\in \{x,y\} \}.    
\end{align*}

It is again easy to exploit McDiarmid's inequality to show that in the random 8-graph $\cH$ we considered earlier
with high probability for all $x,y$ and $j_x,j_y\in \{1,2\}$
$$|\cP_{j_x,j_y}|=   \frac{p^2(1-p)^{10} k^3 n^4}{j_xj_y} \pm O(n^{20/3})$$
and we assume that we have chosen $\cH$ such that this holds now
(indeed, for $w\in V'$, define $b_w=10n^6$ if $w$ is a copy of $x$ or $y$ and $b_w=10n^5$ otherwise; hence $\sum_{w\in V'}b_w^2=O(n^{13})$).
We define an indicator weight function $w_{x,y,j_x,j_y}$ for the pairs in $\cP_{j_x,j_y}$.
Assume for now  that these test functions are $(d, \varepsilon, \cC)$-trackable for all $\eps \in (0,1/4)$.
Adding all these weight functions to the set of weight functions which we give to Theorem~\ref{hypmat}, 
 we obtain a matching $M$ such that 
\begin{equation} \label{eqP2}\left|{M \choose 2} \cap \cP_{j_x,j_y}\right| = w_{x,y,j_x,j_y}(M) \leq (1+d^{-\varepsilon^3})d^{-2}|\cP_{j_x,j_y}| \le (1+n^{-2\delta})\frac{4(1-p)^2 k}{25j_xj_y}\end{equation} for each $x,y,j_x,j_y$.

Now we verify that $w_{x,y,j_x,j_y}$ is $(d, \varepsilon, \cC)$-trackable.
Property \ref{item: trackable size} follows trivially by our estimation of $|\cP_{j_x,j_y}|=\Theta (nd^2)$ above.
To see~\ref{item: trackable degrees}, fix some edge $e$ of $\cH$; as we may assume that this edge is contained in at least some pairs in $\cP_{j_x,j_y}$, 
by symmetry, we may assume $e=(K,i,j)$ 
for some distinct $i,j\in [k]$.
Then the weight of all pairs containing $e$ is $O(d)<nd^{1-1/4}$ as any pair needs to contain at least one vertex in $\{y_i,y_j\}$.
This proves~\ref{item: trackable degrees}.
To see~\ref{item: trackable neighborhood}, we argue as follows.
We select two edges $e=(K_x,\alpha_x, \beta_x),f=(K_y,\alpha_y, \beta_y)$ of $\cH$ so that $\{e,f\}\in \cP_{j_x,j_y}$.
We aim to find an upper bound for the number of triples $g_1,g_2,g_3\in \cH$ such that they form an element of $\cC$ both with $e$ and with $f$.
First we only focus on $e$.
In order to find an upper bound on the number of such triples that form  a conflict with $e$,
we can choose five vertices in $G$ and four colors, resulting in an upper bound of $O(n^9)$.
If the triples $g_1,g_2,g_3$  also form a conflict with $f$,
we can select at most four vertices outside $V(K_x)\cup V(K_y)$ and again at most four colors.
Hence there are at most $O(n^8)<d^{3-1/4}$ such triples.
This proves~\ref{item: trackable neighborhood}.
For future reference, we note here that in order to prove~\ref{item: trackable neighborhood} it suffices to assume that $K_x,K_y$ share at most one vertex.
Property~\ref{item: trackable no conflicts} is vacuously true.

Next we define a set of triples $\cT_{j_x,j_y}$ that extend the pairs in $\cP_{j_x,j_y}$.
Suppose that we have  a pair $\{(K_x,\alpha_x, \beta_x),(K_y,\alpha_y, \beta_y)\}\in \cP_{j_x,j_y}$. 
We add the triple $\{(K_x,\alpha_x, \beta_x),(K_y,\alpha_y, \beta_y), (K, \gamma, \gamma')\}$ to $\cT_{j_x,j_y}$ 
whenever $K$ is edge-disjoint from $K_x,K_y$ and 
$K$ contains vertices $v_x\in V(K_x),v_y\in V(K_y)$ and 
such that $xv_x,yv_y$ have color $i$ in $K_x,K_y$.
The vertex $v_xv_y$ in $\cH$ has degree $d\pm O(n^2)$ as $\cH$ is almost regular and as $\Delta_2(\cH)=O(d/n)$ almost all edges containing $v_xv_y$ avoid $(K_x,\alpha_x, \beta_x),(K_y,\alpha_y, \beta_y)$.
Consequently,  $$|\cT_{j_x,j_y}|=j_xj_y(d \pm O(n^2))|\cP_{j_x,j_y}|.$$
In the same vein as above, we define an indicator weight function $w_{x,y,j_x,j_y}'$ on the triples of $E(\cH)$.
Assume for now that these weight functions are $(d, \varepsilon, \cC)$-trackable for all $\varepsilon \in (0, 1/4)$.
Theorem~\ref{hypmat} gives rise to a matching $M$ such that 
\begin{equation} \label{eqT2}\left|{M \choose 3} \cap \cT_{j_x,j_y}\right|
= w'_{x,y,j_x,j_y}(M) \geq (1-d^{-\varepsilon^3}) d^{-3}|\cT_{j_x,j_y}| \ge (1-n^{-2\delta})\frac{4(1-p)^2 k}{25}
\end{equation}
for each $x,y,j_x,j_y$.

The crucial observation is that each edge $x'y'$ from (V) arises from an edge of $\cH$ in some triple of $\cT_{j_x,j_y}$ that contains a pair from 
$\cP_{j_x,j_y}$. 
Therefore the number of $x'y'$ as in (V) is at most 
$$    \sum_{j_x,j_y\in \{1,3\}} j_xj_y\left|{M\choose 2} \cap \cP_{j_x,j_y}\right|- \left|{M\choose 3} \cap \cT_{j_x,j_y}\right|.$$
By (\ref{eqP2}) and (\ref{eqT2}) we see that this is at most $n^{1-\delta}$, which proves (V).

It remains to check that $w_{x,y,j_x,j_y}'$ is $(d, \varepsilon, \cC)$-trackable for all $\varepsilon \in (0, 1/4)$. 
To see \ref{item: trackable size}, recall that $|\cT_{j_x,j_y}|=\Theta(d|\cP_{j_x,j_y}|)=\Theta(nd^3)$.
To see~\ref{item: trackable degrees}, we follow the same argument as above for~$w_{x,y,j_x,j_y}$
and conclude that there is $O(d^2)<nd^{2-1/4}$ weight on triples containing a particular edge of $\cH$ and 
$O(d)<nd^{1-1/4}$ weight on triples  containing two fixed edges.
As we noted earlier when we checked~\ref{item: trackable neighborhood} for $w_{x,y,j_x,j_y}$,
the number of conflicts that form a conflict with two different edges which have at most one common vertex is at most $d^{3-1/4}$.
This applies to any pair of edges of $\cH$ that belong to a triple with positive weight given by~$w_{x,y,j_x,j_y}'$.
Hence~\ref{item: trackable neighborhood} holds.
Property~\ref{item: trackable no conflicts} is again vacuously true and this completes the proof of Theorem~\ref{thm:56n}.
\qed

To prove~\eqref{eq:56n},
we may proceed very similar as for the proof of~\eqref{eqn:knn} in the previous section and also apply the local lemma
with almost the same collection of bad events.
Here we deviate slightly from the approach in~\cite{BCDP:22} as we aim to use the most basic version of the local lemma which appears easier to us.

First apply Theorem~\ref{thm:56n} (with $2\delta$ playing the role of $\delta$) and obtain a coloring as stated in Theorem~\ref{thm:56n}.
In particular, $\Delta(L)\leq n^{1-2\delta}$ by (IV).
Color the edges of $L$ with a set $P$ of $k=n^{1-\delta}$ new colors independently and uniformly at random.

For any pair $e,f$ of adjacent edges in $L$ and $i\in P$,
we let $A_{e,f,i}$ be the event that both~$e$ and~$f$ receive color $i$.
Clearly, $\bP[A_{e,f,i}]=k^{-2}$.
For a $4$-cycle $D$ in $L$,
we define $B_D$ to be the event that $D$ receives two colors (and this coloring is proper).
Hence $\bP[B_{D}]=(k(k-1))^{-1}\leq 2k^{-2}$.
For a $4$-cycle $D=xyx'y'$ in $K_{n}$ and $i\in P$,
where $xy,x'y'\in E(L)$ and $xy',yx'$ belong to two triangles colored alike by Theorem~\ref{stage1},
we let $C_{D,i}$ be the event that both $xy$, $x'y'$ receive color~$i$.
Clearly, $\bP[C_{D,i}]=k^{-2}$.
Using essentially the same argument verbatim  as in the previous section, we can apply the local lemma to avoid all these events simultaneously.
This shows~\eqref{eq:56n}.\qed

\section{Concluding remarks}
As a further example of the versatility of our method, we can prove the following nonbipartite version of~\eqref{eqn:knn}:
\begin{equation} \label{clique} r(K_n, C_4, 3)=n/2 + O(n^{1-\delta}).\end{equation}
The lower bound is easy: every component of every color class in a $(C_4, 3)$-coloring of $K_n$ is a star or a triangle, hence each color class has at most $n$ edges. So the number of colors is at least ${n \choose 2}/n = (n-1)/2$. For the upper bound, our main tool is the following result that is the nonbipartite analogue of Theorem~\ref{stage1}. 

\begin{theorem} \label{stage1T}
There exists $\delta>0$ such that for all sufficiently large $n$ in terms of $\delta$, 
there exists an edge-coloring of a subgraph of $G \subset K_{n}$ with at most $n/2$ colors  and the following properties:
\begin{enumerate}[label=\rm (\Roman*)]
    \item Every color class consists of vertex-disjoint triangles. 
    \item Given any two colors $i,j$, 
    the 3-graph with vertex set $V(G)$ where each edge is formed by the vertex set of a triangle in color $i$ or $j$ has girth at least 5.
    \item The graph $L=K_n-E(G)$ has maximum degree at most $n^{1-\delta}$. 
    \item For each $xy \in E(K_n)$, the number of  $x'y' \in E(L)$ with $\{x,y\} \cap \{x', y'\}=\emptyset$ for which $xx'$ and $yy'$ receive the same color in $G$ is at most $n^{1-\delta}$.
\end{enumerate}
\end{theorem}

The proof of Theorem~\ref{stage1T} is very similar to the proof of Theorem~\ref{stage1}. One can then quickly use the local lemma to prove (\ref{clique}). 
It is an interesting question to decide if one can improve the error term in this result. 
Roughly speaking, this is equivalent to constructing an essentially resolvable Steiner triple system where the union of any two color classes has girth 5. 
We conjecture that such objects exist.
\begin{conjecture} $r(K_{n}, C_4, 3) = \frac{n}{2} + O(1)$ and $r(K_{n}, C_4, 3) = \frac{n-1}{2}$ for infinitely many $n$.
\end{conjecture}
In a similar vein, we pose another conjecture and problem.
\begin{conjecture} $r(K_{n,n}, C_4, 3) = \frac{2n}{3} + O(1)$ and $r(K_{n,n}, C_4, 3) = \frac{2n}{3}$ for infinitely many $n$.
\end{conjecture}
\begin{problem} Is $r(K_n, K_4, 5) = \frac{5n}{6} + O(1)?$
\end{problem}

\bibliographystyle{amsplain}
\bibliography{references}

\end{document}